\documentclass[12pt]{amsart}
\usepackage{amssymb,latexsym}
\usepackage{enumerate}
\usepackage{hyperref}

\makeatletter
\@namedef{subjclassname@2010}{%
  \textup{2010} Mathematics Subject Classification}
\makeatother

\newtheorem{thm}{Theorem}[section]

\theoremstyle{definition}
\newtheorem{defin}[thm]{Definition}
\newtheorem{exa}[thm]{Example}
\newtheorem{rem}[thm]{Remark}

\newtheorem*{que}{Question}

\DeclareMathOperator\length{length}
\DeclareMathOperator\Tor{Tor}

\begin{document}

\baselineskip=16pt

\title[Serre's Intersection multiplicities]{Computing Serre's Intersection Multiplicities}

\author[Dang Tuan Hiep]{Dang Tuan Hiep}

\address{Dipartimento di Matematica, Universit\`{a} degli Studi di Bari Aldo Moro, Via E. Orabona, 4 - 70125, Bari, Italy}

\email{dang@dm.uniba.it}

\address{{\it Current address:} Fachbereich Mathematik, Universit\"{a}t Kaiserslautern, Postfach 3049, 67653 Kaiserslautern, Germany}

\date{}

\begin{abstract}
The aim of this note is to describe how to compute the intersection multiplicity defined by Jean Pierre Serre. Furthermore, many examples in \cite{Ful} are checked by our implementation in \textsc{Sage} and \textsc{Singular}.
\end{abstract}



\maketitle

\section{Introduction}

There have been many definitions of intersection multiplicities and each has its own range of applications and set of assumptions. Generally the basic idea is to define the order of tangency of two subspaces meeting at one point in such a way that certain natural conditions hold. Serre's definition is purely algebraic and satisfies conditions listed below.

\begin{exa}\label{exam1}
Let consider the intersection of two curves in the plane. One of the curves is the $x$-axis. The other curve is defined by the equation $y = x^3 - x^2$. There are two points of intersection, the origin $(0,0)$ and the point $(1,0)$, the first with multiplicity $2$ and the second with multiplicity $1$. The basic idea behind early algebraic definition of intersection multiplicity is that it should be determined by the dimension of the vector space obtained by dividing the polynomial ring $k[x,y]$ by the ideal generated by the polynomial defining the curves. In this case the polynomials are $y$ and $y - x^3 + x^2$, the quotient $k[x,y]/(y,y-x^3+x^2)$ is isomorphic to $k[x]/(x^3-x^2)$, which has dimension $3$. This number gives the total number of intersections counted with the appropriate multiplicities. In order to obtain the multiplicity at each point, we replace the ring $k[x,y]$ by its localization at that point. For instance, at the point $(0,0)$, if $R$ denotes the localization of the polynomial ring at the ideal $(x,y)$, it is easy to show that the dimension of $R/(y,y-x^3+x^2)$ is $2$. We are able to check these results with \textsc{Sage} \cite{Ste} as follows:
\begin{verbatim}
sage: R.<x,y> = QQ[]; I = R.ideal(y); J = R.ideal(y - x^3 + x^2)
sage: Y = AffineScheme(I); Z = AffineScheme(J)
sage: p1 = (0,0); p2 = (1,0)
sage: Y.intersection_multiplicity(Z, p1)
2
sage: Y.intersection_multiplicity(Z, p2)
1
\end{verbatim}
\end{exa}

\begin{rem}
In order to do the computations using \textsc{Sage} in this note, we need to attach the separate package (\texttt{scheme\_base.py}) using the following command:
\begin{verbatim}
sage: attach scheme_base.py
\end{verbatim}
The source code of this separate package is available at \url{http://uniba-it.academia.edu/HiepDang/Teaching/30154/scheme_base.py}.
\end{rem}

This example suggests that we can define intersection multiplicities in general as follows: we take the local ring at a point, which we denote $R$, take the ideals defining the two subvarieties near the point, say $I$ and $J$, and define the intersection multiplicity as the vector space dimension $\dim_k(R/(I+J))$. We remark that the fact this point is an isolated point of the intersection assures that this dimension is finite. Using this definition we can implement a command to compute the intersection multiplicities as follows:

\begin{enumerate}
\item {\it Affine case:}

\begin{verbatim}
def intersection_multiplicity(self, arg, point):
    local_ring1, trans1 = localize_at_point(self.ring(),point)
    local_ideal1 = trans1(self.ideal())
    local_ring2, trans2 = localize_at_point(arg.ring(),point)
    local_ideal2 = trans2(arg.ideal())
    return AffineScheme(local_ideal1 + local_ideal2).degree()
\end{verbatim}

Note that the command \texttt{localize\_at\_point} is defined below.

\item {\it Projective case:}
\begin{verbatim}
def intersection_multiplicity(self, arg, point):
    R = self.ring()
    i = 0
    for i in range(R.ngens()):
        if point[i] != 0:
            break
    U = self.affine_chart(R.gen(i))
    V = arg.affine_chart(R.gen(i))
    p = list(point)
    pi = p.pop(i)
    new_point = tuple([1/pi*x for x in p])
    return U.intersection_multiplicity(V,new_point)
\end{verbatim}

Note that the command \texttt{affine\_chart} is also defined below.
\end{enumerate}

The problem with this definition is that it lacks some of the properties required by intersection multiplicities; in particular, it does not satisfy B\'{e}zout's Theorem.

\begin{thm}[B\'{e}zout's Theorem, see in \cite{Rob}]
Let $Y$ and $Z$ be closed subschemes of $\mathbb P^n$ such that they are of complementary dimension and intersect in a finite number of points. Then the number of points of intersection counted with multiplicities is the product of the degrees of the $Y$ and $Z$.
\end{thm}

In order to see why the above definition of intersection multiplicity does not satisfy B\'{e}zout's Theorem, we consider the following example.

\begin{exa}[Example 7.1.4 in \cite{Ful}] \label{exa1}
Let $Y$ and $Z$ be subschemes of $\mathbb P^5$ defined by the ideals $I = (xz,xw,yz,yw)$ and $J = (x - z, y - w)$, respectively, where the coordinate ring of $\mathbb P^5$ is $R = \mathbb Q[x,y,z,w,t]$.
\begin{verbatim}
sage: R.<x,y,z,w,t> = QQ[]
sage: I = R.ideal(x*z,x*w,y*z,y*w); J = R.ideal(x - z, y - w)
sage: Y = ProjectiveScheme(I); Z = ProjectiveScheme(J)
sage: point = (0,0,0,0,1)
sage: Y.intersection_multiplicity(Z, point)
3
sage: Y.degree()
2
sage: Z.degree()
1
\end{verbatim}
This means that the intersection of $Y$ and $Z$ is at unique point $(0,0,0,0,1)$ with multiplicity $3$. However the degrees of $Y$ and $Z$ are $2$ and $1$, respectively. Thus B\'{e}zout's Theorem is not satisfied in this case.
\end{exa}

We now give Serre's definition of intersection multiplicities, which does not have this drawback. To improve flexibility, the definition is given in terms of modules $M$ and $N$; the case of subvarieties is the case in which $M = R/I$ and $N = R/J$ as above.

\begin{defin}[\cite{Rob}]
Let $R$ be a regular local ring, and let $M$ and $N$ be finitely generated $R$-modules such that $M \otimes_R N$ is an $R$-module of finite length. Then the intersection multiplicity of $M$ and $N$ is
$$\chi (M,N) = \sum_i (-1)^i \length (\Tor_i^R (M,N)).$$
\end{defin}

This definition requires two conditions that $\Tor_i^R (M,N)$ has finite length for each $i$ and that it be zero for large $i$. The first condition follows from the assumption that $M \otimes_R N$ has finite length. The second condition is a modern version of the Hilbert Syzygy Theorem, that for a regular local ring, the projective dimension of any module is finite.

The intersection multiplicities should have the following properties:
\begin{enumerate}
\item It satisfies B\'{e}zout's Theorem.
\item If $M = R/I$ and $N = R/J$, where $I$ and $J$ are ideals defined by smooth subschemes intersecting transversally, then $\chi(M,N) = 1$.
\end{enumerate}
In the chapter 8 of \cite{Rob}, the author shows that B\'{e}zout's Theorem holds if Serre's definition of intersection multiplicity is used. However, there are other properties that aer not obvious, for example, that the intersection multiplicity is nonnegative. Serre \cite{Ser} conjectured that the following should be true. 

Let $R$ be a regular local ring, and $M$ and $N$ be finitely generated $R$--modules. Suppose that $M\otimes_RN$ is a module of finite length. Then
\begin{enumerate}
\item $\dim(M)+\dim(N) \leq \dim(R)$.
\item $\chi(M,N) \geq 0$.
\item $\chi(M,N)\neq 0$ if and only if $\dim(M)+\dim(N) = \dim(R)$.
\end{enumerate}
It is easy to see that the second and third conditions can be replaced by
\begin{enumerate}
\item (Vanishing) If $\dim(M)+\dim(N) < \dim(R)$, then $\chi(M,N) = 0$.
\item (Positivity) If $\dim(M)+\dim(N) = \dim(R)$, then $\chi(M,N) > 0$.
\end{enumerate}
In fact Serre proved the first one in general and others in many cases. All three were proven in the case of equal characteristic. The details can be found in Serre \cite{Ser}. 

The vanishing conjecture has been proven for arbitrary regular local rings using local Chern characters (Roberts \cite{Rob1, Rob2}) and independently using Adams operations (Gillet and Soul\'{e} \cite{GS}).

The fact $\chi(M,N) \geq 0$ has been proven by Gabber using de Jong's theorem on the existence of regular alterations. Gabber's result can be found in Berthelot \cite{Ber}, and de Jong's theorem can be found in de Jong \cite{deJ}. The positivity conjecture is still open.

\section{How to compute Serre's intersection multiplicities}

Let $Y$ and $Z$ be two subschemes of $\mathbb P^n$. Assume that $Y$ and $Z$ are defined by ideals $I$ and $J$ and that $Y\cap Z$ consists of a finite set of points. For each point $p$ in the intersection, let $\mathfrak m_p$ be the maximal ideal of $R = k[x_0, \ldots, x_n]$ corresponding to $p$, and let $\chi ((R/I)_{\mathfrak m_p}, (R/J)_{\mathfrak m_p})$ over the regular local ring $k[x_0, \ldots, x_n]_{\mathfrak m_p}$ be the intersection multiplicity at $p$.
\begin{que}
How can we compute $\chi ((R/I)_{\mathfrak m_p}, (R/J)_{\mathfrak m_p})$?
\end{que}

\subsection{Affine case}
In the affine case, firstly, we need to return the localization of the polynomial ring at the maximal ideal which corresponding to the intersected point. We also need to return the translation mapping from the polynomial ring to its localization.
\begin{verbatim}
def localize_at_point(ring, point):
    local_ring = PolynomialRing(ring.base_ring(),
                 ring.variable_names(), order = 'ds')
    new_coordinate = [local_ring.gens()[i] + point[i]
                     for i in range(local_ring.ngens())]
    trans = ring.hom(new_coordinate)
    return local_ring, trans
\end{verbatim}
Then we use the interfaces between \textsc{Sage} and \textsc{Singular} \cite{DGPS} to compute
$$\length(\Tor_i^{R_{\mathfrak m_p}}((R/I)_{\mathfrak m_p}, (R/J)_{\mathfrak m_p})),$$
where $R_{\mathfrak m_p}$ denotes the localization of the polynomial ring at the maximal ideal $\mathfrak m_p$ which corresponding to the intersected point $p$.
\begin{verbatim}
def serre_intersection_multiplicity(self, arg, point):
    local_ring, trans = localize_at_point(self.ring(),point)
    I = trans(self.ideal())
    J = trans(arg.ideal())
    from sage.interfaces.singular import singular
    singular.LIB('homolog.lib')
    i = 0
    s = 0
    t = sum(singular.Tor(i, I, J).std().hilb(2).sage())
    while t != 0:
        s = s + ((-1)**i)*t
        i = i + 1
        t = sum(singular.Tor(i, I, J).std().hilb(2).sage())
    return s
\end{verbatim}
\subsection{Projective case}
In the projective case, firstly, we need to return the affine chart of a projective scheme. It should be an affine scheme.
\begin{verbatim}
def affine_chart(self, v):
    ngens = [p.subs({v:1}) for p in self.ideal().gens()]
    L = list(self.ring().gens())
    L.remove(v)
    R = PolynomialRing(self.ring().base_ring(),L)
    return AffineScheme(R.ideal(ngens))
\end{verbatim}
Then we recall the computation in the affine case.
\begin{verbatim}
def serre_intersection_multiplicity(self, arg, point):
    R = self.ring()
    i = 0
    for i in range(R.ngens()):
        if point[i] != 0:
            break
    U = self.affine_chart(R.gen(i))
    V = arg.affine_chart(R.gen(i))
    p = list(point)
    pi = p.pop(i)
    new_point = tuple([1/pi*x for x in p])
    return U.serre_intersection_multiplicity(V, new_point)
\end{verbatim}
Let us consider again the Example \ref{exa1} as follows:
\begin{verbatim}
sage: R.<x,y,z,w,t> = QQ[]
sage: I = R.ideal(x*z,x*w,y*z,y*w); J = R.ideal(x - z, y - w)
sage: Y = ProjectiveScheme(I); Z = ProjectiveScheme(J)
sage: point = (0,0,0,0,1)
sage: Y.serre_intersection_multiplicity(Z, point)
2
\end{verbatim}
Serre's definition shows that the intersection multiplicity of $Y$ and $Z$ at the point $(0,0,0,0,1)$ should be $2$. Thus B\'{e}zout's Theorem is satisfied in this case.

\begin{exa}[Example 7.1.5 in \cite{Ful}]\label{exam2}
In $\mathbb A^4$, let $Y$ be the affine scheme defined by the ideal generated by $x,w$ and let $Z$ be the image of the finite morphism $\varphi$ from $\mathbb A^2$ to $\mathbb A^4$ given by
$$\varphi(s,t) = (s^4,s^3t,st^3,t^4).$$
We use \textsc{Singular} (see in \cite{DL} and \cite{GP} for more details) to compute the ideal defines $Z$ as follows:
\begin{verbatim}
> ring R = 0,(s,t,x,y,z,w),dp;
> ideal J = x-s4,y-s3t,z-st3,w-t4;
> eliminate(eliminate(J,s),t);
_[1]=yz-xw
_[2]=z3-yw2
_[3]=xz2-y2w
_[4]=y3-x2z
\end{verbatim}
Thus $Z$ is an affine scheme defined by the ideal generated by $yz-xw, z^3-yw^2,xz^2-y^2w,y^3-x^2z$. Moreover, it is easy to show that the origin $p = (0,0,0,0)$ is a proper component of the intersection $Y\cap Z$. We use \textsc{Sage} to return the intersection multiplicities as follows:
\begin{verbatim}
sage: R.<x,y,z,w> = QQ[]
sage: I = R.ideal(x,w)
sage: J = R.ideal(y*z-x*w,z^3-y*w^2,x*z^2-y^2*w,y^3-x^2*z)
sage: Y = AffineScheme(I); Z = AffineScheme(J)
sage: p = (0,0,0,0)
sage: Y.intersection_multiplicity(Z,p)
5
sage: Y.serre_intersection_multiplicity(Z,p)
4
\end{verbatim}
\end{exa}

\section{\textsc{Singular} computations}

In this section we present the \textsc{Singular} computations for the above examples. We start by loading all \textsc{Singular} libraries:
\begin{verbatim}
LIB "all.lib";	                     //load all libraries
\end{verbatim}
We write two procedures as follows:
\begin{verbatim}
proc intersection_multiplicity(ideal I, ideal J)
"USAGE:  intersection_multiplicity(I,J); I,J = ideals
RETURN:  the intersection multiplicity of two subvarieties defined by 
the ideals I, J at the origin
"
{
  ideal K = I + J;
  int v = vdim(std(K));
  return (v);
}
proc serre_intersection_multiplicity(ideal I, ideal J)
"USAGE:  serre_intersection_multiplicity(I,J); I,J = ideals
RETURN:  the intersection multiplicity (defined by J. P. Serre) of two 
subvarieties defined by the ideals I, J at the origin
"
{
  int i = 0;
  int s = 0;
  module m = std(Tor(i,I,J));
  int t = sum(hilb(m,2));
  while (t != 0)
  {
    s = s + ((-1)^i)*t;
    i++;
    module m = std(Tor(i,I,J));
    t = sum(hilb(m,2));
  }
  return (s);
}
\end{verbatim}
In order to calculate the intersection multiplicities of two curves in the example \ref{exam1}, we can do as follows:
\begin{verbatim}
ring r = 0, (x,y), ds;
ideal I = y;
ideal J = y-x3+x2;
intersection_multiplicity(I,J);
2               //at the origin (0,0)
ring s = 0, (x,y), ds;
map f = r,x+1,y;
intersection_multiplicity(f(I),f(J));
1               //at the point (1,0)
\end{verbatim}
In order to calculate the intersection multiplicities in the example \ref{exa1}, we can do as follows:
\begin{verbatim}
ring r = 0, (x,y,z,w), ds;
ideal I = xz,xw,yz,yw;
ideal J = x-z,y-w;
intersection_multiplicity(I,J);
3
serre_intersection_multiplicity(I,J);
2
\end{verbatim}
In order to calculate the intersection multiplicities in the example \ref{exam2}, we can do as follows:
\begin{verbatim}
ring r = 0, (x,y,z,w), ds;
ideal I = x,w;
ideal J = yz-xw,z3-yw2,xz2-y2w,y3-x2z;
intersection_multiplicity(I,J);
5
serre_intersection_multiplicity(I,J);
4
\end{verbatim}

\subsection*{Acknowledgements}
The author would like to thank Prof. Dr. Wolfram Decker for many interesting and helpful discussions on the topic of this work. The author also would like to thank PD Dr. Mohamed Barakat and Dr. Bur\c{c}in Er\"{o}cal for the instruction to write \textsc{Sage} code.


\begin{thebibliography}{99}

\bibitem{Ber}
P. Berthelot, Alt\'{e}rations de vari\'{e}t\'{e}s alg\'{e}briques (d'apr\'{e}s A. J. de Jong), \emph{S\'{e}minaire Bourbaki}, (1996), 273--311.

\bibitem{deJ}
A. de Jong, Smoothness, stability, and alterations, \emph{Publ. Math. IHES} {\bf 83} (1996), 51--93.

\bibitem{DL}
W. Decker and C. Lossen, \emph{Computing in Algebraic Geometry: A Quick Start using \textsc{Singular}}, Springer, 2006.

\bibitem{DGPS}
W. Decker, G.-M. Greuel, G. Pfister, H. Sch{\"o}nemann,
\newblock {\textsc{Singular}} {3-1-3} --- {A} computer algebra system for polynomial computations, 2011. Available at \url{http://www.singular.uni-kl.de}.

\bibitem{EGSS}
D. Eisenbud, D. R. Grayson, M. Stillman and B. Sturmfels (Eds.), {\em Computations in Algebraic Geometry with \textsc{Macaulay2}}, Springer, 2002.

\bibitem{GP}
G.-M. Greuel and G. Pfister, \emph{A \textsc{Singular} Introduction to commutative algebra}, second edition, Springer, 2008.

\bibitem{GS}
H. Gillet and C. Soul\'{e}, Intersection theory using Adams operations, \emph{Invent. Math.} {\bf 90} (1987), 243--277.

\bibitem{Ful}
W. Fulton, {\em Intersection Theory}, second edition, Springer, 1997.

\bibitem{Rob}
P. C. Roberts, {\em Multiplicities and Chern classes in Local Algebra}, Cambridge University Press, 1998.

\bibitem{Rob1}
P. C. Roberts, The vanishing of intersection multiplicities of perfect complexes, \emph{Bull. Amer. Math. Soc.} {\bf 13} (1985), 127--130.

\bibitem{Rob2}
P. C. Roberts, Local Chern characters and intersection multiplicities, \emph{Proc.Sympos. Pure. Math.} {\bf 46}(2) (1987), 389--400.

\bibitem{Ser}
J. P. Serre, \emph{Alg\`{e}bre Locale -- multiplicit\'{e}s}, Lecture Notes in Mathematics {\bf 11}, Springer, 1961.

\bibitem{Ste}
W. A. Stein et al. Sage Mathematics Software (Version 4.7), The Sage Development Team, 2011. Available at \url{http://www.sagemath.org}.

\end{thebibliography}
\end{document}